\documentclass{amsart}
\usepackage{amssymb}

\title[On Macdonald and Hall-Littlewood Polynomials]{Combinatorial Formulas for Macdonald and Hall-Littlewood Polynomials of Types $A$ and $C$\\ {\rm{\small{Extended Abstract}}}}

\author{Cristian Lenart}

\address{Department of Mathematics and Statistics, State University of New York at Albany, Albany, NY 12222}
\thanks{The author was partially supported by the National Science Foundation grant  DMS-0701044. The author is grateful to Jim Haglund for helpful discussions.}


\newcommand{\des}{{\rm des}}
\newcommand{\Des}{{\rm Des}}
\newcommand{\Diff}{{\rm Diff}}
\newcommand{\maj}{{\rm maj}}
\newcommand{\inv}{{\rm inv}}
\newcommand{\rt}{{\rm r}}

\newcommand{\arm}{{\rm arm}}
\newcommand{\leg}{{\rm leg}}

\newlength{\cellsize}
\newcommand\tableau[1]{
\vcenter{
\let\\=\cr
\baselineskip=-16000pt
\lineskiplimit=16000pt
\lineskip=0pt
\halign{&\tableaucell{##}\cr#1\crcr}}}

\newcommand{\tableaucell}[1]{{%
\def \arg{#1}\def \void{}%
\ifx \void \arg
\vbox to \cellsize{\vfil \hrule width \cellsize height 0pt}%
\else
\unitlength=\cellsize
\begin{picture}(1,1)
\put(0,0){\makebox(1,1){$#1$}}
\put(0,0){\line(1,0){1}}
\put(0,1){\line(1,0){1}}
\put(0,0){\line(0,1){1}}
\put(1,0){\line(0,1){1}}
\end{picture}%
\fi}}

\newtheorem{theorem}{Theorem}[section]
\newtheorem{proposition}[theorem]{Proposition}

\newtheorem{definition}[theorem]{Definition}

\newtheorem{example}[theorem]{Example}

\newtheorem{remark}[theorem]{Remark}
\newtheorem{remarks}[theorem]{Remarks}

\def\R{\mathbb{R}}
\def\Z{\mathbb{Z}}

\def\Waff{W_{\mathrm{aff}}}

\def\h{\mathfrak{h}}
\def\hR{\mathfrak{h}^*_\mathbb{R}}

\newcommand{\casetwo}[3]{\left\{ \begin{array}{ll} #1 &\mbox{if $#2$} \\ [0.04in] #3 &\mbox{otherwise}\,. \end{array} \right.}

\newcommand{\casetwoex}[4]{\left\{ \begin{array}{ll} #1 &\mbox{if $#2$} \\ #3 &\mbox{if $#4$} \,. \end{array} \right.}
\newcommand{\casetwoexc}[4]{\left\{ \begin{array}{ll} #1 &\mbox{if $#2$} \\ #3 &\mbox{if $#4$} \,, \end{array} \right.}

\begin{document}
\bibliographystyle{plain}

\begin{abstract} 
A breakthrough in the theory of (type $A$) Macdonald polynomials is due to Haglund, Haiman and Loehr, who exhibited a combinatorial formula for these polynomials in terms of fillings of Young diagrams. Recently, Ram and Yip gave a formula for the Macdonald polynomials of arbitrary type in terms of the corresponding affine Weyl group. In this paper, we show that a Haglund-Haiman-Loehr type formula follows naturally from the more general Ram-Yip formula, via compression. Then we extend this approach to the Hall-Littlewood polynomials of type $C$, which are specializations of the corresponding Macdonald polynomials at $q=0$. We note that no analog of the Haglund-Haiman-Loehr formula exists beyond type $A$, so our work is a first step towards finding such a formula.
\end{abstract}

\maketitle

%\vspace{-2mm}

\section{Introduction}
\label{intro}

Macdonald \cite{macsft,macopa} defined a remarkable family of symmetric orthogonal polynomials depending on parameters $q,t$, which bear his name. These polynomials generalize several other symmetric polynomials related to representation theory. For instance, at $q=0$, the Macdonald polynomials specialize to the Hall-Littlewood polynomials (or spherical functions on $p$-adic groups), and they further specialize to the Weyl characters (upon setting $t=0$ as well). There has been considerable interest recently in the combinatorics of Macdonald polynomials. This stems in part from a combinatorial formula for the ones corresponding to type $A$, which is due to Haglund, Haiman, and Loehr \cite{hhlcfm}, and which is in terms of fillings of Young diagrams. This formula uses two statistics on the mentioned fillings, called ``inv'' and ``maj''. The Haglund-Haiman-Loehr formula already found important applications, such as new proofs of the Schur positivity for Macdonald polynomials \cite{asasem,gahaha}. Let us also note that there is a version of the Haglund-Haiman-Loehr formula for the non-symmetric Macdonald polynomials \cite{hhlcfn}, as well as a different formula for these polynomials due to Lascoux \cite{lassam}. 

Schwer \cite{schghl} gave a formula for the Hall-Littlewood polynomials of arbitrary type (cf. also \cite{ramawh}). This formula is in terms of so-called alcove walks, which originate in the work of Gaussent-Littelmann \cite{gallsg} and of the author with Postnikov \cite{lapawg,lapcmc} on discrete counterparts to the Littelmann path model \cite{litlrr,litpro}. Schwer's formula was recently generalized by Ram and Yip to a similar formula for the Macdonald polynomials \cite{raycfm}. The generalization consists in the fact that the latter formula is in terms of alcove walks with both ``positive'' and ``negative'' foldings, whereas in the former only ``positive'' foldings appear.

In this paper, we relate the Ram-Yip formula to the Haglund-Haiman-Loehr formula. More precisely, we show that we can group the terms in the type $A$ version of the Ram-Yip formula into equivalence classes, such that the sum in each class is a term in a new formula, which is similar to the Haglund-Haiman-Loehr one but contains considerably fewer terms. An equivalence class consists of all the terms corresponding to alcove walks that produce the same filling of a Young diagram $\lambda$ (indexing the Macdonald polynomial) via a simple construction. In fact, in this paper we require that the partition $\lambda$ is a regular weight; the general case will be considered elsewhere. 

Our approach has the advantage of deriving the Haglund-Haiman-Loehr statistics ``inv'' and ``maj'' on fillings of Young diagrams in a natural way, from more general concepts. It also has the advantage of being applicable to other root systems, where no analog of the Haglund-Haiman-Loehr formula exists. As a first step in this direction, we derive here a formula in terms of fillings of Young diagrams for the Hall-Littlewood polynomials of type $C$ indexed by a regular weight; we proceed by compressing the type $C$ version of Schwer's formula. We expect that similar formulas exist in types $B$ and $D$.

The structure of this extended abstract is as follows. In Section \ref{sectone} we present our formula of Haglund-Haiman-Loehr type for the Macdonald polynomials of type $A$. In Section \ref{hlpcn} we present our new formula for the Hall-Littlewood polynomials of type $C$ in terms of fillings of Young diagrams. In Section \ref{prelim} we give background information on root systems, alcove walks, the Ram-Yip formula, and Schwer's formula. In Section \ref{compramyip} we specialize the Ram-Yip formula to type $A$ and explain how it compresses to our formula for the corresponding Macdonald polynomials. In Section \ref{compsch} we specialize Schwer's formula to type $C$ and explain how it compresses to our formula for the corresponding Hall-Littlewood polynomials. The full length version of Section \ref{compramyip} is \cite{lencfm}.

\section{A new formula of Haglund-Haiman-Loehr type} \label{sectone}

In this section we present a new formula for the Macdonald polynomials of type $A$ that is similar to the Haglund-Haiman-Loehr one \cite{hhlcfm}. This formula will be derived by compressing the Ram-Yip formula \cite{raycfm}. It also turns out that the new formula has considerably fewer terms even than the Haglund-Haiman-Loehr formula. 

Let us consider a partition with $n-1$ distinct parts $\lambda = (\lambda _{1}> \lambda _{2}> \ldots > \lambda _{n-1}>0)$ for a fixed $n$ (this corresponds to a dominant regular weight for the root system of type $A_{n-1}$). Using standard notation, one defines $n(\lambda):=\sum_{i}(i-1)\lambda_i$. 
We identify $\lambda$ with its Young (or Ferrers) diagram, as usual, and denote by $(i,j)$ the cell in row $i$ and column $j$, where $1\le j\le\lambda_i$. We draw this diagram in ``Japanese style'', that is, we embed it in the third quadrant, as shown below:
\[
\lambda =(4,2) = \tableau{{}&{}&{}&{}\\ &&{}&{}}\;.
\]
For any cell $u=(i,j)$ of $\lambda$ with $j\ne 1$, denote the cell $v=(i,j-1)$ directly to the right of $u$ by $\rt(u)$. 

Two cells $u,v\in \lambda $ are said to \emph{attack} each other if either
\begin{itemize}
\item [(i)] they are in the same column: $u = (i,j)$, $v = (k,j)$; or
\item [(ii)] they are in consecutive columns, with the cell in the left column
strictly above the one in the right column: $u=(i,j)$,
$v=(k,j-1)$, where $i<k$.
\end{itemize}
%The figure below shows the two types of pairs of attacking cells.
%\[
%\quad \tableau{{\bullet}&{}\\ {}&{}\\{\bullet}&{}\\&{}} \, ,\qquad
%\quad  \tableau{{\bullet}&{}\\ {}&{}\\{}&{\bullet}\\&{}}\; .
%\]

\begin{remark} {\rm The main difference in our approach compared to the Haglund-Haiman-Loehr one is in the definition of attacking cells; note that in \cite{hhlcfm} these cells are defined similarly, except that $u=(i,j)$ and $v=(k,j-1)$ with  $i>k$ attack each other.}
\end{remark} 

A \emph{filling} is a function
$\sigma \::\: \lambda \rightarrow [n]:=\{1,\ldots,n\}$ for some $n$, that is, an assignment of values in $[n]$ to the cells of $\lambda$.  As usual, we define the content of a filling $\sigma$ as ${\rm content}(\sigma):=(c_1,\ldots,c_n)$, where $c_i$ is the number of entries $i$ in the filling, i.e., $c_i:=|\sigma^{-1}(i)|$. The monomial $x^{{\rm content}(\sigma)}$ in the variables $x_{1},\ldots,x_n$ is then given by $x^{{\rm content}(\sigma)} := x_1^{c_1}\ldots,x_n^{c_n}$.

\begin{definition}\label{deff} A filling $\sigma\::\:\lambda\rightarrow [n]$ is called \emph{non-attacking} if $\sigma(u)\ne\sigma(v)$ whenever  $u$ and $v$ attack each other. 
Let ${\mathcal T}(\lambda,n)$ denote the set of non-attacking fillings.
\end{definition}

\begin{definition} Given a filling $\sigma$ of $\lambda$, let
\[\Des(\sigma):=\{(i,j)\in\lambda\::\:(i,j+1)\in\lambda\,,\;\;\sigma(i,j)>\sigma(i,j+1)\}\,,\]
%and $\des(\sigma):=|\Des(\sigma)|$. 
\[\Diff(\sigma):=\{(i,j)\in\lambda\::\:(i,j+1)\in\lambda\,,\;\;\sigma(i,j)\ne\sigma(i,j+1)\}\,.\]
\end{definition}

We define a reading order   on the cells of $\lambda$ as  the total order given by considering the columns from right to left (largest to smallest), and by reading each column from top to bottom. Note that this is a different reading order than the usual (French or Japanese) ones.

\begin{definition} An {\em
inversion} of $\sigma $ is a pair $(u,v)$ of attacking cells, where $u$ precedes $v$ in the considered reading order and $\sigma (u)>\sigma (v)$. Let ${\rm Inv}(\sigma)$ denote the set of inversions of $\sigma$.
\end{definition}

Here are two examples of inversions, where $a<b$:
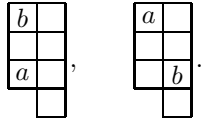
\[
\tableau{{b}&{}\\ {}&{}\\{a}&{}\\&{}} \, ,\qquad 
\tableau{{a}&{}\\ {}&{}\\{}&{b}\\&{}}\,.
\]
%Define
%\begin{equation}\label{e:Inv}
%\Inv (\sigma ) := \{\{u,v \}: \text{$\sigma (u)>\sigma (v)$ is an
%inversion} \}. 
%\end{equation}
%and the inv statistic
%\begin{equation}\label{e:maj-inv}
%\inv (\sigma ) := |\Inv (\sigma )|\,.
%\end{equation}

The \emph{arm} of a cell $u\in\lambda$ is the number of cells strictly to the left of $u$ in the same row; similarly, the \emph{leg} of $u$ is the number of cells strictly below $u$ in the same column.%, as illustrated below.
%\[
%\tableau{{a}&{\bullet}\\ {}&{l}\\{}&{l}\\&{l}} \, ,\qquad \arm(\bullet)=1\,,\;\;\;\leg(\bullet)=3\,.
%\]

\begin{definition} The \emph{maj} and \emph{inv statistics} on fillings $\sigma$ are defined by
\[\maj(\sigma):=\sum_{u\in\Des(\sigma)}\arm(u)\,, \;\;\;\;\inv(\sigma):=|{\rm Inv}(\sigma)|-\sum_{u\in\Des(\sigma)}\leg(u)\,.\]
\end{definition}

We are now ready to state a new combinatorial formula for the Macdonald $P$-polynomials in the variables $X=(x_1,\ldots,x_n)$. 

\begin{theorem}\label{hlqthm}
Given a partition $\lambda$ with $n-1$ distinct parts, we have
\begin{equation}\label{hlqform}
P_{\lambda}(X;q,t) = \sum_{\sigma\in{\mathcal T}(\lambda,n)}
 t^{n(\lambda) - \inv (\sigma)}q^{\maj(\sigma)}\left(\prod_{u\in\Diff(\sigma)}\frac{1-t}{1-q^{\arm(u)}t^{\leg(u)+1}}\right) x^{{\rm content}(\sigma) } \,.
\end{equation}
\end{theorem}

\section{Hall-Littlewood polynomials of type $C_n$}\label{hlpcn}
%%%%%%%%%%%%%%%n>=2.

In this section we present a new formula for the Hall-Littlewood polynomials of type $C$ in terms of fillings of Young diagrams. This formula will be derived by compressing Schwer's formula \cite{schghl} (cf. also \cite{ramawh}).

Let $\lambda=(\lambda_1>\ldots>\lambda_n>0)$ be a partition with $n$ distinct parts for a fixed $n\ge 2$  (this corresponds to a dominant regular weight for the root system of type $C_n$). Consider the shape $\widehat{\lambda}$ obtained from $\lambda$ by replacing each column of height $k$ with $k$ or $2k-1$ (adjacent) copies of it, depending on the given column being the first one or not. We are representing a filling $\sigma$ of $\widehat{\lambda}$ as a concatenation of columns $C_{ij}$ and $C'_{ik}$, where $i=1,\ldots,\lambda_1$, while for a given $i$ we have $j=1,\ldots,\lambda_i'$ if $i>1$, $j=1$ if $i=1$, and $k=2,\ldots,\lambda_i'$; the columns $C_{ij}$ and $C_{ik}'$ have height $\lambda_i'$. The diagram $\widehat{\lambda}$ is represented in ``Japanese style'', like in the previous section, i.e., the heights of columns increase from left to right; more precisely, we let
\[\sigma={\mathcal C}^{\lambda_1}\ldots {\mathcal C}^{1}\,,\;\;\;\;\mbox{where}\;\;{\mathcal C}^{i}:=\casetwoex{C_{i2}'\ldots C_{i,\lambda_i'}'C_{i1}\ldots C_{i,\lambda_i'}}{i>1}{C_{i2}'\ldots C_{i,\lambda_i'}'C_{i1}}{i=1}\]
Note that the leftmost column is  $C_{\lambda_1,1}$, and the rightmost column is $C_{11}$. For an example, we refer to Section \ref{compsch}.

Essentially, the above description says that the column to the right of $C_{ij}$ is $C_{i,j+1}$, whereas the column to the right of $C_{ik}'$ is $C_{i,k+1}'$. Here we are assuming that the mentioned columns exist, up to the following conventions:
\begin{equation}\label{conv}C_{i,\lambda_i'+1}=\casetwoexc{C_{i-1,2}'}{i>1\:\mbox{ and }\:\lambda_{i-1}'>1}{C_{i-1,1}}{i>1\:\mbox{ and }\:\lambda_{i-1}'=1}\;\;\;\;\;C_{i,\lambda_i'+1}'=C_{i1}\,.\end{equation}

We consider the alphabet $[\overline{n}]:=\{1<\ldots<n<\overline{n}<\overline{n-1}<\ldots<\overline{1}\}$, where the barred entries are viewed as negatives, so that $-\overline{\imath}=i$. Next, we consider the set ${\mathcal T}(\widehat{\lambda},\overline{n})$ of fillings of $\widehat{\lambda}$ with entries in $[\overline{n}]$ which satisfy the following conditions:
\begin{enumerate}
\item the rows are weakly decreasing from left to right;
\item no column contains two entries $a,b$ with $a=\pm b$;
\item any two adjacent columns are related as indicated below.
\end{enumerate}

In order to explain the mentioned relation between adjacent columns, we consider  right actions of type $C$ reflections on columns (see Section \ref{compsch}). For instance, $C(a,\overline{b})$ is the column obtained from $C$ by transposing the entries in positions $a,b$ and by changing their signs. Let us first explain the passage from some column $C_{ij}$ to $C_{i,j+1}$. There exist positions $1\le r_1<\ldots <r_p<j$ (possibly $p=0$) such that $C_{i,j+1}$ differs from $D=C_{ij}(r_1,\overline{\jmath})\ldots (r_p,\overline{\jmath})$ only in position $j$, while  $C_{i,j+1}(j)\not\in \{\pm D(r)\::\:r\in[\lambda_i']\setminus\{j\}\}$ and $C_{i,j+1}(j)\le D(j)$. To include the case $j=\lambda_i'$ in this description, just replace $C_{i,j+1}$ everywhere by $C_{i,j+1}[1,\lambda_i']$ and use the conventions (\ref{conv}). Let us now explain the passage from some column $C_{ik}'$ to $C_{i,k+1}'$. There exist positions $1\le r_1<\ldots <r_p<k$ (possibly $p=0$) such that $C_{i,k+1}'=C_{ik}'(r_1,\overline{k})\ldots (r_p,\overline{k})$. This description includes the case $k=\lambda_i'$, based on the conventions (\ref{conv}).

Let us now define the content of a filling. For this purpose, we first associate with a filling $\sigma$ a compressed version of it, namely the filling $\overline{\sigma}$ of the partition $2\lambda$. This is defined as follows:
\begin{equation}\label{redfilling}
\overline{\sigma}=\overline{{\mathcal C}}^{\lambda_1}\ldots \overline{{\mathcal C}}^{1}\,,\;\;\;\;\mbox{where $\:\overline{{\mathcal C}}^{i}:=C_{i2}'C_{i1}$}\,,\end{equation}
where the conventions (\ref{conv}) are used again. Now define ${\rm content}(\sigma)=(m_1,\ldots,m_n)$, where $m_i$ is half the difference between the number of occurences of the entries $i$ and $\overline{\imath}$ in $\overline{\sigma}$. 

We now define two statistics on fillings that will be used in our compressed formula for Hall-Littlewood polynomials. Intervals refer to the discrete set $[\overline{n}]$. 
Let 
\[\sigma_{ab}:=\casetwo{1}{a,b\ge\overline{n}}{0}\]
Given  a sequence of integers $w$, we write $w[i,j]$ for the subsequence $w(i)w(i+1)\ldots w(j)$. We use the notation $N_{ab}(w)$ for the number of entries $w(i)$ in the interval $(a,b)$. 

Given two columns $D,C$ of the same height $d$ such that $D\ge C$ in the componentwise order, we will define two statistics $N(D,C)$ and $\des(D,C)$ in some special cases, as specified below.

\emph{Case} 0. If $D=C$, then $N(D,C):=0$ and $\des(D,C):=0$. 

\emph{Case} 1. Assume that $C=D(r,\overline{\jmath})$ with $r<j$. Let $a:=D(r)$ and $b:=D(j)$. In this case, we set
\[N(D,C):=N_{\overline{b}a}(D[r+1,j-1])+|(\overline{b},a)\setminus\{\pm D(i)\::\:i=1,\ldots,j\}|+\sigma_{ab}\,,\;\;\;\des(D,C):=1\,.\]
%\[N(D,C):=N_{\overline{b}a}-N_{\overline{b}a}(D[1,r-1])-N_{\overline{a}b}(D[1,r,j-1])+\sigma_{ab}\,,\;\;\;\des(D,C):=1\,.\]

\emph{Case} 2. Assume that $C=D(r_1,\overline{\jmath})\ldots (r_p,\overline{\jmath})$ where $1\le r_1<\ldots <r_p<j$. Let $D_i:=D(r_1,\overline{\jmath})\ldots (r_i,\overline{\jmath})$ for $i=0,\ldots,p$, so that $D_0=D$ and $D_p=C$. We define
\[N(D,C):=\sum_{i=1}^p N(D_{i-1},D_i)\,,\;\;\;\des(D,C):=p\,.\]

\emph{Case} 3. Assume that $C$ differs from $D':=D(r_1,\overline{\jmath})\ldots (r_p,\overline{\jmath})$ with $1\le r_1<\ldots <r_p<j$ (possibly $p=0$) only in position $j$, while  $C(j)\not\in \{\pm D'(r)\::\:r\in[d]\setminus\{j\}\}$ and $C(j)\le D'(j)$. We define
\[N(D,C):=N(D,D')+N_{C(j),D'(j)}(D[j+1,d])\,,\;\;\;\des(D,C):=p+\delta_{C(j),D'(j)}\,,\]
where $\delta_{a,b}$ is the Kronecker delta.

If the height of $C$ is larger than the height $d$ of $D$ (necessarily by 1), and $N(D,C[1,d])$ can be computed as above, we let $N(D,C):=N(D,C[1,d])$ and $\des(D,C):=\des(D,C[1,d])$. Given a filling $\sigma$ with columns $C_m,\ldots ,C_1$, we set
\[N(\sigma):=\sum_{i=1}^{m-1}N(C_{i+1},C_i)+\inv(C_1)\,,\]
assuming that all the terms $N(\,\cdot\,,\,\cdot\,)$ on the right-hand side are of the types described above; here $\inv(C_1)$ denotes the number of (ordinary) inversions in $C_1$, that is, the number of pairs $i<j$ of positions in $C_1$ with $C_1(i)>C_1(j)$. Furthermore, in the mentioned case, we also set 
 \[\des(\sigma):=\sum_{i=1}^{m-1}\des(C_{i+1},C_i)\,.\]
%Note that $\des(\sigma)$ essentially counts the descents in the rows of $\sigma$. 

We can now state our new formula for the Hall-Littlewood polynomials of type $C$. We refer to Remarks \ref{specialfill} for more comments on this formula. 

\begin{theorem}\label{newformc}
Given a partition $\lambda$ with $n$ distinct parts, the Hall-Littlewood polynomial $P_\lambda(X;t)$ is given by
\begin{equation}\label{newhhlc}P_\lambda(X;t)=\sum_{\sigma\in{\mathcal T}(\widehat{\lambda},\overline{n})} t^{N(\sigma)}\,(1-t)^{\des(\sigma)}\,x^{\rm content(\sigma)}\,.\end{equation}
\end{theorem}

\section{Alcove walks and Macdonald polynomials}\label{prelim}

\subsection{Root systems}\label{rootsyst}

We recall some background information on finite root systems and affine Weyl groups.
Let $\mathfrak{g}$ be a complex semisimple Lie algebra, and $\h$ a Cartan subalgebra, whose rank is $r$.
Let $\Phi\subset \h^*$ be the 
corresponding irreducible \emph{root system}, $\hR\subset \h^*$ the real span of the roots, and $\Phi^+\subset \Phi$ the set of positive roots. Let $\rho:=\frac{1}{2}(\sum_{\alpha\in\Phi^+}\alpha)$. 
Let $\alpha_1,\ldots,\alpha_r\in\Phi^+$ be the corresponding 
\emph{simple roots}.
We denote by $\langle\,\cdot\,,\,\cdot\,\rangle$ the non-degenerate scalar product on $\hR$ induced by
the Killing form.  
Given a root $\alpha$, we consider the corresponding \emph{coroot\/} $\alpha^\vee := 2\alpha/\langle\alpha,\alpha\rangle$ and reflection $s_\alpha$.  %The collection of coroots $\Phi^\vee:=\{\alpha^\vee \::\: \alpha\in\Phi\}$ forms the \emph{dual root system.}

Let $W$ be the corresponding  \emph{Weyl group\/}, whose Coxeter generators are denoted, as usual, by $s_i:=s_{\alpha_i}$. The length function on $W$ is denoted by $\ell(\,\cdot\,)$. The \emph{Bruhat order} on $W$ is given by its covers $w\lessdot ws_\beta$, where $\beta\in\Phi^+$, and $\ell(ws_\beta)=\ell(w)+1$.

The \emph{weight lattice\/} $\Lambda$ is given by $\Lambda:=\{\lambda\in \hR \::\: \langle\lambda,\alpha^\vee\rangle\in\Z\}$ for any $\alpha\in\Phi$. The weight lattice $\Lambda$ is generated by the 
\emph{fundamental weights\/}
$\omega_1,\ldots,\omega_r$, which form the dual basis to the 
basis of simple coroots, i.e., $\langle\omega_i,\alpha_j^\vee\rangle=\delta_{ij}$.
The set $\Lambda^+$ of \emph{dominant weights\/} is given by $\Lambda^+:=\{\lambda\in\Lambda \::\: \langle\lambda,\alpha^\vee\rangle\geq 0\}$ for any $\alpha\in\Phi^+$.
%The subgroup of $W$ stabilizing a weight $\lambda$ is denoted by $W_\lambda$, and the set of minimum coset representatives in $W/W_\lambda$ by $W^\lambda$. 
Let $\Z[\Lambda]$ be the group algebra of the weight lattice $\Lambda$, which  has
a $\Z$-basis of formal exponents $\{x^\lambda \::\: \lambda\in\Lambda\}$ with
multiplication $x^\lambda\cdot x^\mu := x^{\lambda+\mu}$.

Given  $\alpha\in\Phi$ and $k\in\Z$, we denote by $s_{\alpha,k}$ the reflection in the affine hyperplane
\begin{equation}
H_{\alpha,k} := \{\lambda\in \hR \::\: \langle\lambda,\alpha^\vee\rangle=k\}.
\label{eqhyp}
\end{equation}
These reflections generate the \emph{affine Weyl group\/} $\Waff$ for the \emph{dual root system} 
$\Phi^\vee:=\{\alpha^\vee \::\: \alpha\in\Phi\}$. 
The hyperplanes $H_{\alpha,k}$ divide the real vector space $\hR$ into open
regions, called \emph{alcoves.} 
The \emph{fundamental alcove\/} $A^\circ$ is given by 
$$
A^\circ :=\{\lambda\in \hR \::\: 0<\langle\lambda,\alpha^\vee\rangle<1 \textrm{ for all }
\alpha\in\Phi^+\}.
$$

\subsection{Alcove walks}\label{alcovewalks}

We say that two alcoves $A$ and $B$ are \emph{adjacent} 
if they are distinct and have a common wall.  
Given a pair of adjacent alcoves $A\ne B$ (i.e., having a common wall), we write 
$A\stackrel{\beta}\longrightarrow B$ if the common wall 
is of the form $H_{\beta,k}$ and the root $\beta\in\Phi$ points 
in the direction from $A$ to $B$.  

\begin{definition} {\rm \cite{lapawg}}
An \emph{alcove path\/} is a sequence of alcoves
 such that any two consecutive ones are adjacent. 
We say that an alcove path $(A_0,A_1,\ldots,A_m)$ is \emph{reduced\/} if $m$ is the minimal 
length of all alcove paths from $A_0$ to $A_m$.
\end{definition}

We need the following generalization of alcove paths.

\begin{definition}\label{defalcwalk} An \emph{alcove walk} is a sequence 
$\Omega=(A_0,F_1,A_1, F_2, \ldots , F_m, A_m, F_{\infty})$ 
such that $A_0,\ldots,$ $A_m$ are alcoves; %$F_{0}$ is a vertex of the first alcove $A_0$; 
$F_i$ is a codimension one common face of the alcoves $A_{i-1}$ and $A_i$,
for $i=1,\ldots,m$; and 
$F_{\infty}$ is a vertex of the last alcove $A_m$. The weight $F_\infty$ is called the \emph{weight} of the alcove walk, and is denoted by $\mu(\Omega)$. 
%The weight $\mu$ is called the \emph{weight\/} of the alcove walk and is denoted by $\mu(\Omega)$.
\end{definition}
 
The \emph{folding operator} $\phi_i$ is the operator which acts on an alcove walk by leaving its initial segment from $A_0$ to $A_{i-1}$ intact and by reflecting the remaining tail in the affine hyperplane containing the face $F_i$. In other words, we define
$$\phi_i(\Omega):=(A_0, F_1, A_1, \ldots, A_{i-1}, F_i'=F_i,  A_{i}', F_{i+1}', A_{i+1}', \ldots,  A_m', F_{\infty}')\,,$$
where $A_j' := \rho_i(A_j)$ for $j\in\{i,\ldots,m\}$, $F_j':=\rho_i(F_j)$ for $j\in\{i,\ldots,m\}\cup\{\infty\}$, and $\rho_i$ is the affine reflection in the hyperplane containing $F_i$. Note that any two folding operators commute. An index $j$ such that $A_{j-1}=A_j$ is called a \emph{folding position} of $\Omega$. Let $\mbox{fp}(\Omega):=\{ j_1<\ldots< j_s\}$ be the set of folding positions of $\Omega$. If this set is empty, $\Omega$ is called \emph{unfolded}. Given this data, we define the operator ``unfold'', producing an unfolded alcove walk, by
\[\mbox{unfold}(\Omega)=\phi_{j_1}\ldots \phi_{j_s} (\Omega)\,.\]

\begin{definition} 
A folding position $j$ of the alcove walk $\Omega=(A_0,F_1,A_1, F_2, \ldots , F_m, A_m, F_{\infty})$ is called a \emph{positive folding} if the alcove $A_{j-1}=A_j$ lies on the positive side of the affine hyperplane containing the face $F_j$. Otherwise, the folding position is called a \emph{negative folding}.
\end{definition}

Let $\tau_\lambda\in\Waff$ denote the translation by $\lambda$. Recall the bijection $A\mapsto v_A$ between alcoves and affine Weyl group elements given by $v_A(A^\circ)=A$. We now fix a dominant weight $\lambda$ and a reduced alcove path $\Pi:=(A_0,A_1,\ldots,A_m)$ from $A^\circ=A_0$ to the alcove $A_m$ corresponding to the minimal element in the coset $\tau_\lambda W$ under the mentioned bijection. Assume that we have
\begin{equation}\label{lambdachain}A_0\stackrel{\beta_1}\longrightarrow A_1\stackrel{\beta_2}\longrightarrow \ldots
\stackrel{\beta_m}\longrightarrow A_{m}\,,\end{equation}
where $\Gamma:=(\beta_1,\ldots,\beta_m)$ is a sequence of positive roots. This sequence, which determines the alcove path, is called a \emph{$\lambda$-chain} (of roots). 

%\begin{remark} $\lambda$-chains were defined in \cite{lapcmc,lapawg} based on alcove paths from $A^\circ$ to $A^\circ-\lambda$. Two equivalent definitions of such $\lambda$-chains (in terms of reduced words in affine Weyl groups, and an interlacing condition) can be found in \cite{lapawg}[Definition 5.4] and \cite{lapcmc}[Definition 4.1 and Proposition 4.4]. Hence, the $\lambda$-chains considered in this paper  are obtained by reversing the ones in the mentioned papers and by removing a certain segment at the end. The reason for this removal is that the alcove paths here do not end at the alcove $A^\circ + \lambda$, but at the minimum length representative in its orbit.
%\end{remark}

 We also let $r_i:=s_{\beta_i}$, and let $\widehat{r}_i$ be the affine reflection in the common wall of $A_{i-1}$ and $A_i$, for $i=1,\ldots,m$; in other words, $\widehat{r}_i:=s_{\beta_i,l_i}$, where $l_i:=|\{j\le i\::\: \beta_j = \beta_i\}|$ is the cardinality of the corresponding set. Given $J=\{j_1<\ldots<j_s\}\subseteq[m]:=\{1,\ldots,m\}$, we define the Weyl group element $\phi(J)$ and the weight $\mu(J)$ by
\begin{equation}\label{defphimu}\phi(J):={r}_{j_1}\ldots {r}_{j_s}\,,\;\;\;\;\;\mu(J):=\widehat{r}_{j_1}\ldots \widehat{r}_{j_s}(\lambda)\,.\end{equation}

\subsection{The Ram-Yip formula for Macdonald polynomials}  Given $w\in W$ and the alcove path $\Pi$ considered above, we define the alcove path 
\[w(\Pi):=(w(A_0),w(A_1),\ldots,w(A_m))\,.\]
Consider the set of alcove paths ${\mathcal P}(\Gamma):=\{w(\Pi)\::\:w\in W\}$. We identify any $w(\Pi)$ with the obvious unfolded alcove walk of weight $\mu(w(\Pi)):=w(\lambda)$. Let us now consider the set of alcove walks
\[{\mathcal F}(\Gamma):=\{\,\mbox{alcove walks $\Omega$}\::\:\mbox{unfold}(\Omega)\in{\mathcal P}(\Gamma)\}\,.\]
We can encode an alcove walk $\Omega$ in ${\mathcal F}(\Gamma)$ by the pair $(w,J)$ in $W\times 2^{[m]}$, where 
\[\mbox{fp}(\Omega)=J\;\;\;\;\mbox{and}\;\;\;\;\mbox{unfold}(\Omega)=w(\Pi)\,.\]
Clearly, we can recover $\Omega$ from $(w,J)$ with $J=\{j_1<\ldots<j_s\}$ by $\Omega=\phi_{j_1}\ldots \phi_{j_s} (w(\Pi))$. 
We call a pair $(w,J)$ a \emph{folding pair}, and, for simplicity, we denote the set $W\times 2^{[m]}$ of such pairs by ${\mathcal F}(\Gamma)$ as well. Given a folding pair $(w,J)$, the corresponding positive and negative foldings (viewed as a partition of $J$) are denoted by $J^+$ and $J^-$.   

\begin{proposition}\label{admpairs} {\rm (1)} Consider a folding pair $(w,J)$ with $J=\{j_1<\ldots<j_s\}$. We have $j_i\in J^+$ if and only if $wr_{j_1}\ldots r_{j_{i-1}}>wr_{j_1}\ldots r_{j_{i-1}}r_{j_{i}}$. {\rm (2)} If $\Omega\mapsto (w,J)$, then $\mu(\Omega)=w(\mu(J))$.
\end{proposition}

We call the sequence $w,\,wr_{j_1},\,\ldots,\,wr_{j_1}\ldots r_{j_{s}}=w\phi(J)$ the Bruhat chain associated to $(w,J)$. 

We now restate the Ram-Yip formula \cite{raycfm} for the Macdonald polynomials $P_\lambda(X;q,t)$ in terms of folding pairs. From now on we assume that the weight $\lambda$ is regular (and dominant), i.e., $\langle\lambda,\alpha^\vee\rangle>0$ for all positive roots $\alpha$. 

\begin{theorem}{\rm \cite{raycfm}} \label{hlpthm} Given a dominant regular weight $\lambda$, we have (based on the notation in Section {\rm \ref{alcovewalks}})
\begin{equation}\label{hlpform}\!\!\!\!\!\!\!\!\!\!\!\!\!\!\!\!\!\!\!\!\!\!\!\!\!\!\!\!\!\!\!\!\!\!\!\!\!\!\!\!\!\!\!\!\!\!\!\!\!\!\!\!\!\!\!\!\!\!\!\!\!\!\!\!\!\!\!\!\!\!\!\!\!\!\!\!\!\!\!\!\!\!\!\!\!\!\!\!\!\!\!\!\!\!\!\!\!\!\!\!\!\!\!\!\!\!\!\!\!\!\!\!\!\!\!\!\!\!\!\!\!\!\!\!\!\!\!\!\!\!\!\!\!\!\!\!\!\!\!\!\!\!\!\!\!\!\!\!\!\!\!\!\!\!\!\!\!\!\!\!\!\!\!\!\!\!\!\!\!\!\!\!\!\!\!\!\!\!\!\!\!\!\!\!\!\!\!\!\!\!\!\!\!\!\!\!\!\!\!\!\!\!\!\!\!\!\!P_{\lambda}(X;q,t)=\end{equation}
\[=\sum_{(w,J)\in{\mathcal F}(\Gamma)}t^{\frac{1}{2}(\ell(w)-\ell(w\phi(J))-|J|)}\,(1-t)^{|J|}\left(\prod_{j\in J^+}\frac{1}{1-q^{l_j}t^{\langle\rho,\beta_j^\vee\rangle}}\right)\left(\prod_{j\in J^-}\frac{q^{l_j}t^{\langle\rho,\beta_j^\vee\rangle}}{1-q^{l_j}t^{\langle\rho,\beta_j^\vee\rangle}}\right)x^{w(\mu(J))}\,.\]
\end{theorem}

\subsection{Schwer's formula for Hall-Littlewood polynomials}\label{subsschwer} Let us now consider a reduced alcove path from $A^\circ$ to $A^\circ+\lambda$. The associated chain of roots $\Gamma$, defined as in (\ref{lambdachain}), will be called an \emph{extended $\lambda$-chain}. All the previous definitions can be adapted to this setup. Let ${\mathcal F}_+(\Gamma)$ consist of the folding pairs $(w,J)$ with $J_-=\emptyset$, which will be called \emph{positive folding pairs}. 

\begin{theorem}{\rm \cite{ramawh,schghl}} \label{hlpthmh} Given a dominant regular weight $\lambda$, the Hall-Littlewood polynomial $P_\lambda(X;t)$ is given by
\begin{equation}\label{hlpformh}P_{\lambda}(X;t)=\sum_{(w,J)\in{\mathcal F}_+(\Gamma)}t^{\frac{1}{2}(\ell(w)+\ell(w\phi(J))-|J|)}\,(1-t)^{|J|}\,x^{w(\mu(J))}\,.\end{equation}
\end{theorem}

\section{Compressing the Ram-Yip formula in type $A_{n-1}$}\label{compramyip}

We now restrict ourselves to the root system of type $A_{n-1}$, fow which the Weyl group $W$ is the symmetric group $S_n$. Permutations $w\in S_n$ are written in one-line notation $w=w(1)\ldots w(n)$. 
We can identify the space $\h_\R^*$ with the quotient space 
$V:=\R^n/\R(1,\ldots,1)$,
where $\R(1,\ldots,1)$ denotes the subspace in $\R^n$ spanned 
by the vector $(1,\ldots,1)$.  
The action of the symmetric group $S_n$ on $V$ is obtained 
from the (left) $S_n$-action on $\R^n$ by permutation of coordinates.
Let $\varepsilon_1,\ldots,\varepsilon_n\in V$ 
be the images of the coordinate vectors in $\R^n$.
The root system $\Phi$ can be represented as 
$\Phi=\{\alpha_{ij}:=\varepsilon_i-\varepsilon_j \::\: i\ne j,\ 1\leq i,j\leq n\}$.
The simple roots are $\alpha_i=\alpha_{i,i+1}$, 
for $i=1,\ldots,n-1$.
The fundamental weights are $\omega_i = \varepsilon_1+\ldots +\varepsilon_i$, 
for $i=1,\ldots,n-1$. 
The weight lattice is $\Lambda=\Z^n/\Z(1,\ldots,1)$. A dominant weight $\lambda=\lambda_1\varepsilon_1+\ldots+\lambda_{n-1}\varepsilon_{n-1}$ is identified with the partition $(\lambda _{1}\geq \lambda _{2}\geq \ldots \geq \lambda _{n-1}\geq\lambda_n=0)$ of length at most $n-1$. We fix such a partition $\lambda$ for the remainder of this section.%paper, and assume that the corresponding weight is regular, i.e., $(\lambda _{1}> \lambda _{2}> \ldots > \lambda _{n-1} >\lambda_n=0)$. %This means that $\lambda'$ has $m_i>0$ parts equal to $i$ for each $i=1,\ldots,n-1$. 

For simplicity, we use the same notation $(i,j)$ with $i<j$ for the root $\alpha_{ij}$ and the reflection $s_{\alpha_{ij}}$, which is the transposition of $i$ and $j$.  Consider the following chain of roots, denoted by $\Gamma(k)$:
\begin{equation}\label{omegakchain}\begin{array}{lllll}
(&\!\!\!\!(k,n),&(k,n-1),&\ldots,&(k,k+1)\,,\\
&\!\!\!\!(k-1,n),&(k-1,n-1),&\ldots,&(k-1,k+1)\,,\\
&&&\ldots\\
&\!\!\!\!(1,n),&(1,n-1),&\ldots,&(1,k+1)\,\,)\,.
\end{array}\end{equation}
Denote by $\Gamma'(k)$ the chain of roots obtained by removing the root $(i,k+1)$ at the end of each row. Now define a chain $\Gamma$ as a concatenation $\Gamma:=\Gamma_{\lambda_1}\ldots\Gamma_2$, where 
\[\Gamma_j:=\casetwo{\Gamma'(\lambda'_j)}{\mbox{$j=\min\;\{i\::\:\lambda_i'=\lambda_j'\}$}}{\Gamma(\lambda'_j)}\]
%Based on the interlacing condition in \cite{lapcmc}[Definition 4.1 and Proposition 4.4], 
It is not hard to verify that $\Gamma$ is a $\lambda$-chain in the sense discussed in Section \ref{alcovewalks}. The $\lambda$-chain $\Gamma$ is fixed for the remainder of this section. Thus, we can replace the notation ${\mathcal F}(\Gamma)$ with ${\mathcal F}(\lambda)$.

\begin{example}\label{ex21} {\rm Consider $n=4$ and $\lambda =(4,3,1,0)$, for which we have the following $\lambda$-chain (the underlined pairs are only relevant in Example \ref{ex21c} below):
\begin{equation}\label{exlchain}\Gamma=\Gamma_4\Gamma_3\Gamma_2=(\underline{(1,4)},(1,3)\:|\:(2,4),\underline{(2,3)},(1,4),\underline{(1,3)}\:|\:\underline{(2,4)},(1,4))\,.\end{equation}
%We represent the Young diagram of $\lambda$ inside a broken $4\times 4$ rectangle, as shown below. In this way, a transpositions $(i,j)$ in $\Gamma$ can be viewed as swapping entries in the two parts of each column (in rows $i$ and $j$, where the row numbers are also indicated below). 
%\[
% \begin{array}{l} \tableau{{1}&{1}&{1}&{1}\\ &{2}&{2}&{2}\\&&&{3}}\\ \\
%\tableau{{2}\\ {3}&{3}&{3}\\ {4}&{4}&{4}&{4}} \end{array}
%\]
%
}
\end{example}

Given the $\lambda$-chain $\Gamma$ above, in Section \ref{alcovewalks} we considered subsets $J=\{ j_1<\ldots< j_s\}$ of $[m]$, where $m$ is the length of the $\lambda$-chain. Instead of $J$, it is now convenient to use the subsequence of $\Gamma$ indexed by the positions in $J$. This is viewed as a concatenation with distinguished factors $T=T_{\lambda_1}\ldots T_2$ induced by the factorization of $\Gamma$ as $\Gamma_{\lambda_1}\ldots\Gamma_2$. The partition $J=J^+\sqcup J^-$ induces partitions $T=T^+\sqcup T^-$ and $T_j=T_j^+\sqcup T_j^-$. All the notions defined in terms of $J$ are now redefined in terms of $T$. As such, from now on we will write  $\phi(T)$, $\mu(T)$, and $|T|$, the latter being the size of $T$. If $(w,J)$ is a folding pair, we will use the same name for the corresponding pair $(w,T)$. We will use the notation ${\mathcal F}(\Gamma)$ and ${\mathcal F}(\lambda)$ accordingly.  We denote by $wT_{\lambda_1}\ldots T_{j}$ the permutation obtained from $w$ via right multiplication by the transpositions in $T_{\lambda_1},\ldots, T_{j}$, considered from left to right. This agrees with the above convention of using pairs to denote both roots and the corresponding reflections. As such, $\phi(J)$ in (\ref{defphimu}) can now be written simply $T$. 

\begin{example}\label{ex21c}{\rm We continue Example \ref{ex21}, by picking the folding pair $(w,J)$ with $w=2341\in S_4$ and $J=\{1,4,6,7\}$ (see the underlined positions in (\ref{exlchain})). Thus, we have
\[T=T_4T_3T_2=((1,4)\:|\:(2,3),(1,3)\:|\:(2,4))\,.\]
Note that $J^+=\{1,7\}$ and $J^-=\{4,6\}$. Indeed, we have the following Bruhat chain associated to $(w,T)$, where the transposed entries are shown in bold (we represent permutations as broken columns):%, as discussed in Example \ref{ex21}):
\[w=\begin{array}{l}\tableau{{{\mathbf 2}}} \\ \\ \tableau{{3}\\{4}\\{{\mathbf 1}}} \end{array} \:>\:\begin{array}{l}\tableau{{1}} \\ \\ \tableau{{3}\\{4}\\{2}} \end{array}\:|\: \begin{array}{l}\tableau{{{ 1}}\\{{\mathbf 3}}} \\ \\ \tableau{{{\mathbf 4}}\\{2}} \end{array}\:<\: \begin{array}{l}\tableau{{{\mathbf 1}}\\{{4}}} \\ \\ \tableau{{{\mathbf 3}}\\{2}}\end{array}\:<\: \begin{array}{l}\tableau{{3}\\{4}}\\ \\ \tableau{{1}\\{2}}\end{array}\:|\: 
\begin{array}{l}\tableau{{3}\\{{\mathbf 4}}}\\ \\ \tableau{{1}\\{{\mathbf 2}}}\end{array} 
\:>\: \begin{array}{l} \tableau{{3}\\{2}} \\ \\ \tableau{{1}\\{4}}
\end{array} \:|\: \begin{array}{l} \tableau{{3}\\{2}\\{1}} \\ \\ \tableau{{4}}
\end{array} \,.\]
}
\end{example}

Given a folding pair $(w,T)$, we consider the permutations
\[\pi_j=\pi_j(w,T):=wT_{\lambda_1}T_{\lambda_1-1}\ldots T_{j+1}\,,\]
for $j=1,\ldots,\lambda_1$. In particular, $\pi_{\lambda_1}=w$. 

\begin{definition}\label{deffill}
The \emph{filling map} is the map $f$ from folding pairs $(w,T)$ to fillings $\sigma=f(w,T)$ of the shape $\lambda$, defined by $\sigma(i,j):=\pi_j(i)$. 
%In other words, the $j$-th column of the filling $\sigma$ (from right to left) consists of the first $\lambda_j'$ entries of the permutation $\pi_j$.
\end{definition}

\begin{example} {\rm Given $(w,T)$ as in Example \ref{ex21c}, we have
\[f(w,T)=\tableau{{2}&{1}&{3}&{3}\\&{3}&{4}&{2}\\&&&{1}}\,.\]}
\end{example}

From now on, we assume that the partition $\lambda$  corresponds to a regular weight, i.e., $(\lambda_1>\ldots>\lambda_{n-1}>0)$. We will now describe the way in which the formula (\ref{hlqform}) of Haglund-Haiman-Loehr type can be obtained by compressing the Ram-Yip formula (\ref{hlpform}). Thus, Theorem \ref{hlqthm} becomes a corollary of Theorem \ref{mainthm} below. We start by rewriting the Ram-Yip formula (\ref{hlpform}) in the type $A$ setup, as follows: 
\[
 P_{\lambda}(X;q,t)=\sum_{(w,T)\in{\mathcal F}(\Gamma)}t^{\frac{1}{2}(\ell(w)-\ell(wT)-|T|)}\,(1-t)^{|T|}\!\!\left(\prod_{j,(i,k)\in T_j^+}\frac{1}{1-q^{\arm(-j+1,-i)}t^{k-i}}\right)\times\]
\[\;\;\;\;\;\;\;\;\;\;\;\;\;\;\times\left(\prod_{j,(i,k)\in T_j^-}\frac{q^{\arm(-j+1,-i)}t^{k-i}}{1-q^{\arm(-j+1,-i)}t^{k-i}}\right)x^{w(\mu(T))}\,.\]

\begin{theorem}\label{mainthm}
We have $f(\mathcal{F}(\lambda))={\mathcal T}(\lambda,n)$. Given any $\sigma \in{\mathcal T}(\lambda,n)$ and any $(w,T)\in f^{-1}(\sigma)$, we have ${\rm content}(f(w,T))=w(\mu(T))$. Furthermore, the following compression formula holds for any $\sigma \in{\mathcal T}(\lambda,n)$:
\[\sum_{(w,T)\in f^{-1}(\sigma)}\!\!\!\!\!\!\!\!\!t^{\frac{1}{2}(\ell(w)-\ell(wT)-|T|)}\,(1-t)^{|T|}\left(\prod_{j,(i,k)\in T_j^+}\frac{1}{1-q^{\arm(i,j-1)}t^{k-i}}\right)\times\]
\[ \times\left(\prod_{j,(i,k)\in T_j^-}\frac{q^{\arm(i,j-1)}t^{k-i}}{1-q^{\arm(i,j-1)}t^{k-i}}\right)=t^{n(\lambda) - \inv (\sigma)}q^{\maj(\sigma)}\left(\prod_{u\in\Diff(\sigma)}\frac{1-t}{1-q^{\arm(u)}t^{\leg(u)+1}}\right) \,.\]
\end{theorem}

\section{Compressing Schwer's formula in type $C_n$}\label{compsch} We now restrict ourselves to the root system of type $C_n$. 
We can identify the space $\h_\R^*$ with  
$V:=\R^n$, the coordinate vectors being $\varepsilon_1,\ldots,\varepsilon_n$.  
The root system $\Phi$ can be represented as 
$\Phi=\{\pm\varepsilon_i\pm\varepsilon_j \::\:  1\leq i<j\leq n\}\cup\{\pm 2\varepsilon_i\::\: 1\leq i\leq n\}$. 
The simple roots are $\alpha_i=\varepsilon_i-\varepsilon_{i+1}$, 
for $i=1,\ldots,n-1$ and $\alpha_n=2\varepsilon_n$. 
The fundamental weights are $\omega_i = \varepsilon_1+\ldots +\varepsilon_i$, 
for $i=1,\ldots,n$. 
The weight lattice is $\Lambda=\Z^n$. A dominant weight $\lambda=\lambda_1\varepsilon_1+\ldots+\lambda_n\varepsilon_n$ is identified with the partition $(\lambda _{1}\geq \lambda _{2}\geq \ldots \geq \lambda_n\geq 0)$ of length at most $n$. We fix such a partition $\lambda$ for the remainder of this section.

The corresponding Weyl group $W$ is the group of signed permutations $B_n$. Such permutations are bijections $w$ from $[\overline{n}]:=\{1<\ldots<n<\overline{n}<\overline{n-1}<\ldots<\overline{1}\}$ to $[\overline{n}]$ satisfying $w(\overline{\imath})=\overline{w(i)}$. We use the window notation $w=w(1)\ldots w(n)$. The group $B_n$ acts on $V$ as usual, by permuting the coordinate vectors and by changing their signs. 

For simplicity, we use the same notation $(i,j)$ with $1\le i<j\le n$ for the positive root $\varepsilon_i-\varepsilon_j$ and the corresponding reflection, which, in the window notation, is the transposition of entries in positions $i$ and $j$. Similarly, we denote by $(i,\overline{\jmath})$, again for $i<j$, the positive root $\varepsilon_i+\varepsilon_j$ and the corresponding reflection; in the window notation, the latter is the transposition of entries in positions $i$ and $j$ followed by the sign change of those entries. Finally, we denote by $(i,\overline{\imath})$ the positive root $2\varepsilon_i$ and the corresponding reflection, which is the sign change in position $i$.

Let 
\[\Gamma(k)=\Gamma_2'\ldots\Gamma_k'\Gamma_1(k)\ldots\Gamma_k(k)\,,\]
where
\[
\;\:\Gamma_j':=((1,\overline{\jmath}),(2,\overline{\jmath}),\ldots,(j-1,\overline{\jmath}))\,,\]
\[\begin{array}{llllll}\Gamma_j(k):=(\!\!\!\!\!&(1,\overline{\jmath}),&(2,\overline{\jmath}),&\ldots,&(j-1,\overline{\jmath}),\\
&(j,\overline{k+1}),&(j,\overline{k+2}),&\ldots,&(j,\overline{n}),&(j,\overline{\jmath}),\\ 
&(j,n),&(j,n-1),&\ldots,&(j,k+1)\,)\,.\end{array}\]
It is not hard to see that $\Gamma(k)$ is an extended $\omega_k$-chain, in the sense discussed in Section \ref{subsschwer}.
Hence, we can construct an extended $\lambda$-chain as a concatenation $\Gamma:=\Gamma^{\lambda_1}\ldots\Gamma^1$, where 
\begin{equation}\label{fact}\Gamma^i=\Gamma(\lambda'_i)=\Gamma_{i2}'\ldots\Gamma_{i,\lambda_i'}'\Gamma_{i1}\ldots\Gamma_{i,\lambda_i'}\,,\;\;\;\;\mbox{and }\:\Gamma_{ij}=\Gamma_j(\lambda_i')\,,\;\:\Gamma_{ij}'=\Gamma_j'\,.\end{equation}
This extended $\lambda$-chain is fixed for the remainder of this section. Thus, we can replace the notation ${\mathcal F}_+(\Gamma)$ with ${\mathcal F}_+(\lambda)$.

\begin{example}\label{ex21h} {\rm Consider $n=3$ and $\lambda =(3,2,1)$, for which we have the extended $\lambda$-chain below. The factorization of $\Gamma$ into subchains is indicated by vertical bars, while the double vertical bars separate the subchains corresponding to different columns. The underlined pairs are only relevant in Example \ref{ex21ch} below.
\begin{equation}\label{exlchainh}\!\!\!\!\!\!\!\!\!\!\!\!\!\!\!\!\!\!\!\!\!\!\!\!\!\!\!\!\!\!\!\!\!\!\!\!\!\!\!\!\!\!\!\!\!\!\!\!\!\!\!\!\!\!\!\!\!\!\!\!\!\!\!\!\!\!\!\!\!\!\!\!\!\!\!\!\!\!\!\!\!\!\!\!\!\!\!\!\!\!\!\!\!\!\!\!\Gamma=\Gamma_{31}\:||\:\Gamma_{22}'\Gamma_{21}\Gamma_{22}\:||\:\Gamma_{12}'\Gamma_{13}'\Gamma_{11}\Gamma_{12}\Gamma_{13}=\end{equation}\[=((1,\overline{2}), \underline{(1,\overline{3})},(1,\overline{1}),(1,3),(1,2)\:||\:\underline{(1,\overline{2})}\:|\:(1,\overline{3}),(1,\overline{1}),(1,3)\:|\:(1,\overline{2}),(2,\overline{3}),\underline{(2,\overline{2})},\underline{(2,3)}\:|| \]\[\!\!\!\!\!\!\!\!\!\!\!\!\!\!\!\!\!\!\!\!\!\!\!\!\!\!\!\!\!\!\!\!\!\!\!\!\!\!\!\!\!\!\!\!\!\!\!\!\!\!\!\!\!\!\!(1,\overline{2})\:|\:(1,\overline{3}),(2,\overline{3})\:|\:(1,\overline{1})\:|\:(1,\overline{2}), (2,\overline{2})\:|\:(1,\overline{3}),(2,\overline{3}), (3,\overline{3}))\,.\]
%We represent the Young diagram of $\lambda$ inside a broken $3\times 3$ rectangle, as below. In this way, a reflection in $\Gamma$ can be viewed as swapping entries and/or  changing signs in the two parts of each column, or only the top part. 
%\[
% \begin{array}{l} \tableau{{1}&{1}&{1}\\ &{2}&{2}\\&&{3}}\\ \\
%\tableau{{2}\\ {3}&{3}} \end{array}
%\]
}
\end{example}

Given the extended $\lambda$-chain $\Gamma$ above, in Section \ref{alcovewalks} we considered subsets $J=\{ j_1<\ldots< j_s\}$ of $[m]$, where  $m$ is the length of $\Gamma$. Instead of $J$, it is now convenient to use the subsequence of $\Gamma$ indexed by the positions in $J$. This is viewed as a concatenation with distinguished factors $T_{ij}$ and $T_{ik}'$ induced by the factorization (\ref{fact}) of $\Gamma$. All the notions defined in terms of $J$ are now redefined in terms of $T$. As such, from now on we will write  $\phi(T)$, $\mu(T)$, and $|T|$, the latter being the size of $T$. If $(w,J)$ is positive folding pair, we will use the same name for the corresponding pair $(w,T)$.  %We will use the notation ${\mathcal F}_+(\Gamma)$ and ${\mathcal F}_+(\lambda)$ accordingly.  
We denote by $wT_{\lambda_1,1}\ldots T_{ij}$ and $wT_{\lambda_1,1}\ldots T_{ik}'$ the permutations obtained from $w$ via right multiplication by the reflections in $T_{\lambda_1,1},\ldots, T_{ij}$ and $T_{\lambda_1,1},\ldots, T_{ik}'$, considered from left to right. This agrees with the above convention of using pairs to denote both roots and the corresponding reflections. As such, $\phi(J)$ in (\ref{defphimu}) can now be written simply $T$.

\begin{example}\label{ex21ch}{\rm We continue Example \ref{ex21h}, by picking the positive folding pair $(w,J)$ with $w=\overline{1}\,\overline{2}\,\overline{3}\in B_3$ and $J=\{2,6,12,13\}$ (see the underlined positions in (\ref{exlchainh})). Thus, we have 
\[T=T_{31}\:||\:T_{22}'T_{21}T_{22}\:||\:T_{12}'T_{13}'T_{11}T_{12}T_{13}=((1,\overline{3})\:||\:(1,\overline{2})\:|\;\;\;|\:(2, \overline{2}),(2,3)\:||\;\;\;|\;\;\;|\;\;\;|\;\;\;|\;\;\;)\,.\]
We have the following decreasing Bruhat chain associated to $(w,T)$, where the modified entries are shown in bold (we represent signed permutations as broken columns, as in Example \ref{ex21c}, and we display the splitting of the chain into subchains induced by the above splitting of $T$):
\[w=
\begin{array}{l}\tableau{{\mathbf{\overline{1}}}}\\ \\ \tableau{{\overline{2}}\\{\mathbf{\overline{3}}}} \end{array}\:>\:
\begin{array}{l} \tableau{{3}}\\ \\ \tableau{{\overline{2}}\\{1}} \end{array}\:||\:
\begin{array}{l} \tableau{{\mathbf{3}}\\{\mathbf{\overline{2}}}}\\ \\ \tableau{{1}} \end{array}\:>\:
\begin{array}{l} \tableau{{2}\\{\overline{3}}}\\ \\ \tableau{{1}} \end{array}\:|\:
\begin{array}{l} \tableau{{2}\\{\overline{3}}}\\ \\ \tableau{{1}} \end{array}\:|\:
\begin{array}{l} \tableau{{2}\\{\mathbf{\overline{3}}}}\\ \\ \tableau{{1}} \end{array}\:>\:
\begin{array}{l} \tableau{{2}\\{\mathbf{3}}}\\ \\ \tableau{{\mathbf{1}}} \end{array}\:>\:
\begin{array}{l} \tableau{{2}\\{1}}\\ \\ \tableau{{3}} \end{array}\:||\:
\begin{array}{l} \tableau{{2}\\{1}\\{3}}\\ \\ \end{array} \:|\:
\begin{array}{l} \tableau{{2}\\{1}\\{3}}\\ \\ \end{array} \:|\:
\begin{array}{l} \tableau{{2}\\{1}\\{3}}\\ \\ \end{array} \:|\:
\begin{array}{l} \tableau{{2}\\{1}\\{3}}\\ \\ \end{array} \:|\:
\begin{array}{l} \tableau{{2}\\{1}\\{3}}\\ \\ \end{array}
\,.
\]}
\end{example}

Given a positive folding pair $(w,T)$, with $T$ split into factors $T_{ij}$ and $T_{ik}'$ as above, we consider the signed permutations
\[\pi_{ij}=\pi_{ij}(w,T):=wT_{\lambda_1,1}\ldots T_{i,j-1}\,,\;\;\;\;\;\pi_{ik}'=\pi_{ik}'(w,T):=wT_{\lambda_1,1}\ldots T_{i,k-1}'\,;\]
when undefined, $T_{i,j-1}$ and $T_{i,k-1}'$ are given by conventions similar to (\ref{conv}), based on the corresponding factorization (\ref{fact}) of the extended $\lambda$-chain $\Gamma$. In particular, $\pi_{\lambda_1,1}=w$. 

Let us now recall the notation in Section \ref{hlpcn}.

\begin{definition}\label{deffillh}
The \emph{filling map} is the map $\widehat{f}$ from positive folding pairs $(w,T)$ to fillings $\sigma=\widehat{f}(w,T)$ of the shape $\widehat{\lambda}$, defined by $C_{ij}=\pi_{ij}[1,\lambda_i']$ and $C_{ik}'=\pi_{ik}'[1,\lambda_i']$.
\end{definition}

\begin{example} {\rm Given $(w,T)$ as in Example \ref{ex21ch}, we have
\[\widehat{f}(w,T)=\tableau{{\overline{1}}&{3}&{2}&{2}&{2}&{2}&{2}\\&{\overline{2}}&{\overline{3}}&{\overline{3}}&{1}&{1}&{1}\\&&&&{3}&{3}&{3}}\,.\]}
\end{example}

From now on, we assume that the partition $\lambda$  corresponds to a regular weight, i.e., $(\lambda_1>\ldots>\lambda_{n}>0)$. We will now describe the way in which the formula (\ref{newhhlc}) can be obtained by compressing Schwer's formula (\ref{hlpformh}). Thus, Theorem \ref{newformc} becomes a corollary of the theorem below. 

\begin{theorem}\label{ccompress}
We have $\widehat{f}(\mathcal{F_+}(\lambda))={\mathcal T}(\widehat{\lambda},\overline{n})$. Given any $\sigma \in{\mathcal T}((\widehat{\lambda},\overline{n})$ and $(w,T)\in \widehat{f}^{-1}(\sigma)$, we have $w(\mu(T))={\rm content}(\widehat{f}(w,T))$. Furthermore, the following compression formula holds for any $\sigma \in{\mathcal T}(\widehat{\lambda},\overline{n})$:
\begin{equation}\label{fcomp}\sum_{(w,T)\in \widehat{f}^{-1}(\sigma)}t^{\frac{1}{2}(\ell(w)+\ell(w\phi(T))-|T|)}\,(1-t)^{|T|}=t^{N(\sigma)}\,(1-t)^{\des(\sigma)}\,.\end{equation}
\end{theorem}

\begin{remarks}\label{specialfill} {\rm The ``doubled'' versions of the Kashiwara-Nakashima tableaux \cite{kancgr} of shape $\lambda$, which index the basis elements of the irreducible representation of $\mathfrak{sp}_{2n}$ of highest weight $\lambda$, are among the fillings $\overline{\sigma}$ (see (\ref{redfilling})), for $\sigma\in{\mathcal T}(\widehat{\lambda},\overline{n})$. Indeed, it was proved in \cite{laaapm} that for each Kashiwara-Nakashima tableau there is a unique positive folding pair $(w,T)$ whose associated Bruhat chain is saturated and ends at the identity, such that the compressed version $\overline{\widehat{f}(w,T)}$ of $\widehat{f}(w,T)$ is the ``doubled'' version of the given tableau. 

(2) Consider a filling $\tau$ of $\lambda$ which satisfies the following conditions: (i) the rows are weakly decreasing from left to right; (ii) two entries $a,b$ with $a=\pm b$ cannot appear in the same column, or in two consecutive columns in positions $(i,j)$ and $(k,j-1)$ with $i>k$. Let $\tau^2$ be the filling of $2\lambda$ obtained by doubling each column of $\tau$. It is not hard to see that there is a unique filling $\sigma$ in ${\mathcal T}(\widehat{\lambda},\overline{n})$ such that its compressed version $\overline{\sigma}$ coincides with $\tau^2$. For such fillings, we can describe the statistic $N(\sigma)$ in (\ref{fcomp}) in terms of $\tau$, in a similar way to the statistic $\inv$ of Haglund-Haiman-Loehr type in (\ref{hlqform}).}
\end{remarks}

\end{document}